\documentclass[10pt]{amsart}
\usepackage{amsmath,amsfonts,amsthm}
\usepackage{graphicx}
\usepackage{epsfig}
\def\ds{\displaystyle}

\def\Z{{\mathbb Z}}

\def\N{{\mathbb N}}
\def\ml{{\mathcal L}}
\def\md{{\mathcal D}}
\def\eps{\varepsilon}
\def\d{\delta}
\def\t{\tilde}
\author[N. Enriquez]{Nathana\"el ENRIQUEZ}
\address{Laboratoire de Probabilit\'es de Paris 6, 4
place Jussieu,
 75252 Paris cedex 05}
\email{enriquez@ccr.jussieu.fr}

\title[Fractional Brownian motion.]{A simple construction of the 
fractional Brownian motion.}

\begin{document}

\begin{abstract}
   In this work we introduce correlated random walks on $\Z$. When picking
suitably at random the coefficient of correlation, and taking the average over a large
number of walks, we obtain a discrete Gaussian process, whose scaling limit is the
fractional Brownian motion. We have to use two radically different models for both cases
${1\over2}\leq H<1$ and $0<H<{1\over2}$. This result provides an algorithm for
the simulation of the fractional Brownian motion, which appears to be quite efficient.
\end{abstract}
\maketitle
\newtheorem{defi}{Definition}
\newtheorem{theo}{Theorem}
\newtheorem{prop}{Proposition}
\newtheorem{cor}{Corollary}
\newtheorem{lem}{Lemma}
\newtheorem{rem}{Remark}
\section{Introduction}
The fractional Brownian motion appears to be a very natural object,
for its three fundamental and characteristic features of being a continuous Gaussian
process,  self-similar and with stationary increments. By self-similar, we mean that there
exists a real number $H\in]0,1[$, such that the finite-dimensional distributions of $\{
T^{-H}B_H(Tt),\, t\geq0\}$ do not depend on $T$. The parameter $H$ is called the Hurst
parameter or the index of self-similarity. Under these conditions, it is not hard to check
that the covariance function must have the following form:
$$E[B_H(s)B_H(t)]={1\over2}(s^{2H}+t^{2H}-|s-t|^{2H}).$$ 
A major consequence of this fact is
that if we introduce the sequence of increments $G_H(j)=B_H(j)-B_H(j-1), j=1,2,...$, also
called fractional Gaussian noise, we note that they are strongly correlated (for $H\neq
{1\over2}$). More precisely, 
$$ E[G_H(j)G_H(j+k)]\mathop\sim_{k\to\infty} H(2H-1)k^{2H-2}\qquad(1)$$
We denote the two radically different behaviours on both sides of ${1\over2}$:
for $H<1/2$ the increments are all negatively correlated, which corresponds to a chaotic
behavior, whereas for $H>1/2$ the positive correlation between the increments corresponds
to a more disciplined behaviour. We refer to Chapter 7 of Samorodnitsky and Taqqu
\cite{ST} for more information about fractional Brownian motion.

There is a well-known representation introduced by Mandelbrot and Van Ness \cite{MVN} of the
fractional Brownian motion as an integral of a kernel function with respect to the usual
Brownian motion. The method of approximation that uses the discretization of this integral
leads to an algorithm obliging to store a lot of data in memory, and to deal with non
smooth functions in the case $0<H<{1\over2}$.  This inconvenience is nicely discussed by
Carmona and Coutin in \cite{CC}, and they propose a way to reduce it. 

The construction we propose is based on correlated random walks: it consists in
discrete processes such that the law of each move is ruled by the value of the previous
move. We refer to \cite {Enr} and the references therein for more information about these
processes. We note that the decay of correlation for such processes is exponential, but
that a mixture of these processes, adapted to the value of the index of self-similarity,
leads to walks whose correlations satisfy (1).  The superposition of a large number of such
walks yields a discrete Gaussian process whose correlations fulfill the conditions of Taqqu
\cite {Taq}, so that its scaling limit is the fractional Brownian motion. 
We have to distinguish both cases ${1\over2}\leq H<1$ and
$0<H<{1\over2}$. In the first case, the equivalence relation (1) is the only
condition to check. The case $0<H<{1\over2}$ is more delicate: a compensation relation
between the correlations has to be satisfied simultanously. That leads us to introduce two
different  types of correlated random walks for both cases.

This construction reminds the constructions using  renewal processes of
Mandelbrot \cite{M} and Taqqu and Levy \cite {TL},\cite{LT}. It reminds even more the
construction  coming from traffic modeling by Taqqu, Willinger and Sherman
\cite{TWS}. To our knowledge, all these methods are restricted to the case
${1\over2}\leq H<1$. In that case, we shall discuss the differences between our construction
and \cite {TWS}. Finally, we want to mention that random media has already been used
by Kesten and Spitzer \cite{KS} to get convergences towards certain self-similar processes
with index of self-similarity bigger than ${1\over2}$.

In the last part, we discuss the speed of the convergence stated in the
previous parts. We indicate how to deduce a quite fast algorithm that does not ask to keep a
lot of datas in memory. 
\section{The case $1/2\leq H< 1$}

\subsection{The correlated random walk}
We introduce first our basic tool: the correlated random walk with persistence $p$.
It is a process evolving on $\Z$ by jumps of +1 or -1, whose probability of making the same
jump as the previous one is $p$:

\begin{defi}
For any $p\in[0,1]$, the correlated random walk $X^p$ with persistence $p$ is a $\Z$-valued
discrete process, such that :

\_ $X_0^p=0$, $P(X_1^p=-1)=1/2$, $P(X_1^p=1)=1/2$. 

\_ $\forall n\geq1, \, \eps_n^p:=X_{n}^p-X_{n-1}^p$ equals 1 or -1 a.s.

\_ $\forall n\geq1, P(\eps_{n+1}^p=\eps_n^p|\sigma(X_k^p, 0\leq k\leq n))=p.$

\end{defi}
This process is not Markovian, but if we define a state as the position of the process on
$\Z$, coupled with the sign of its last jump, we have to deal with a
Markov process on $\Z\times\{-1,1\}$. In fact, this process consists in alternate 
falls and rises of i.i.d. geometric length with parameter $p$.

We can compute the correlations between two steps
distant from $n$:

\begin{prop}
 $\forall m\geq1, n\geq0,$
$E[\eps_{m}^p\eps_{m+n}^p]=(2p-1)^n  $
\end{prop}

Proof: $\forall n\geq1$,
 
$E[\eps_{n+1}^p|\sigma(X_k, 0\leq k \leq n)]=$

$E[\eps_{n+1}^p1_{\eps_n^p=1}|\sigma(X_k, 0\leq k \leq
n)]+E[\eps_{n+1}^p1_{\eps_n^p=-1}|\sigma(X_k, 0\leq k
\leq n)]= $

$(2p-1) 1_{\eps_n^p=1}-(2p-1)1_{\eps_n^p=-1}=(2p-1)\eps_n^p  $.

Consequently, conditioning by $\sigma(X_k, 0\leq k \leq m+n)$, we get

$\forall m\geq1, n\geq0,$ $E[\eps_{m+n+1}^p\eps_m^p]=(2p-1)E[\eps_{m+n}^p\eps_m^p]$. 

The result is then obtained by
recurrence. \qed

We now introduce an extra randomness in the persistence.
We first denote by $P^p$ the law of
$X^p$ for a given $p$. Now, considering a probability measure $\mu$ on $[0,1]$, 
we will call $P^\mu$, the annealed law of the correlated walk associated to
$\mu$, i.e. the measure on $\Z^\N$ defined by $P^\mu:=\int_0^1 P^p d\mu(p) $. 

Remark: unlike the situation in \cite{Enr}, the persistence does not depend on the
level. Only one coin toss, according to $\mu$, decides for the whole environment.

Let $X^\mu$ be a process of law $P^\mu$.
Let us now introduce the notation $\eps^\mu_n:=
X_{n}^\mu-X_{n-1}^\mu$. From Proposition 1 we get the straightforward result: 
\begin{prop}
$\forall m\geq1, n\geq0,$
$E[\eps_{m}^\mu\eps_{m+n}^\mu]=\int_0^1(2p-1)^nd\mu(p). $
\end{prop}

\subsection{Statement and proof of the result}
The goal now is to introduce a probability measure $\mu$ leading to the same equivalent as
(1), mentioned in the introduction, so that by taking the average over a large number of
trajectories, we approximate a discrete Gaussian process having the same properties as in
\cite{Taq}, whose scaling limit is the fractional Brownian motion:

\begin{theo}
Let $H\in]1/2,1[$.
\par\noindent
Denote by $\mu^H$ the probability on $[{1\over2},1]$ with density
$(1-H)2^{3-2H}(1-p)^{1-2H}$.

Let $(X^{\mu^H,i})_{i\geq1}$  be a sequence of independent processes of law $P^{\mu^H}$,
$$ \ml^\md\lim_{N\to\infty}\ml\lim_{M\to\infty} c_H
{X^{\mu^H,1}_{[Nt]}+...+X^{\mu^H,M}_{[Nt]}\over N^H\sqrt{M}} = B_H(t)$$
with $c_H=\sqrt{H(2H-1)\over\Gamma(3-2H)}$,

$\ml$ means convergence in the sense of the finite-dimensional distributions, and
$\ml^\md$ means convergence in the sense of the weak convergence in the Skorohod topology
on $D[0,1]$, the space of cadlag functions on [0,1].
\end{theo}

Proof: 
The central limit theorem implies that $ \ml\lim_{M\to\infty} 
{X^{\mu^H,1}_{k}+...+X^{\mu^H,M}_{k}\over \sqrt{M}}$ is a discrete centered Gaussian process
$(Y^H_k)_{k\geq1}$, with stationary increments $G^H_k:=Y^H_{k+1}-Y^H_k$ with $E[G^H_k]=0$,
$E[(G^H_k)^2]=1$ and 
$$\forall i,n\geq0,\quad r(n):=E[G^H_i G^H_{i+n}]=
(2-2H)2^{2-2H}\int_{1\over2}^1(2u-1)^n(1-u)^{1-2H}du$$ 
$\begin{array}{rl}
r(n)&=\ds(2-2H) \int_0^1v^n(1-v)^{1-2H}dv\\
&=\ds(2-2H){\Gamma(n+1)\Gamma(2-2H)\over\Gamma(n+3-2H)}+O({1\over n})\\
&\ds\mathop{\sim}_{n\to\infty} \Gamma(3-2H){1\over n^{2-2H}}
={1\over c_H^2}{H(2H-1)\over n^{2-2H}}\end{array}$

So that,
 
$\begin{array}{rl}E[c_H^2(G^H_1+...+G^H_N)^2]&=\ds c_H^2\sum_{i=1}^N\sum_{j=1}^N r(|i-j|)\\
&=
\ds c_H^2(r(0)+\sum_{i=1}^{N-1}[r(0)+2\sum_{k=1}^i r(k)])\\
&\ds\mathop{\sim}_{n\to\infty} N^{2H}\end{array}$

(the last step consists simply in two successive comparisons between sums and integrals).
A direct application of \cite{Taq} (lemma 5.1) allows to conclude.
\qed

We can give also an analog statement for $H=1/2$:

\begin{theo} 

Denote by $\mu^{1\over2}$ the uniform probability on $[{1\over2},1]$.

Let $(X^{\mu^{1\over2},i})_{i\geq1}$  be a sequence of independent processes of law
$P^{\mu^{1\over2}}$,
$$ \ml^\md\lim_{N\to\infty}\ml\lim_{M\to\infty}c_{1\over2}
{X^{\mu^{1\over2},1}_{[Nt]}+...+X^{\mu^{1\over2},M}_{[Nt]}\over \sqrt{N\log N}\sqrt{M}} =
B(t)$$
where $B$ is the classical Brownian motion, and $c_{1\over2}={1\over\sqrt{2}}$.
\end{theo}
Proof: The scheme is the same as in Theorem 1. The difference here is that 
$r(n)=2\int_{1\over2}^1(2u-1)^ndu={1\over (n+1)}$.

So that,
$r(0)+\sum_{i=1}^{N-1}[r(0)+2\sum_{k=1}^i r(k)]\ds\mathop{\sim}_{n\to\infty}2 N\log
N$. 

We conclude again, applying \cite{Taq} (lemma 5.1). \qed

Remark: The order of the limits in both theorems is of big importance: the limit in the
reverse order would bring 0, as far as for any fixed $p$, a correlated random walk
satisfies a central limit theorem with normalization $\sqrt N=o(\sqrt{N\log N})$ (see
\cite{Enr}). 

We want now to compare previous theorems with the result of \cite{TWS}: in \cite{TWS}, the
limit theorem consists also in taking the scaling limit of the average over a large number
(tending to infinity) of i.i.d. copies of processes which are a succession of falls and
rises. In the case of \cite{TWS}, falls and rises are all independent with infinite
variance. 

In our setting, even if the lengths of falls and rises have finite variance under each
$P^p$  (they are geometric), the laws under $P^{\mu^H}$ of the falls and the rises (which
are all the same) have infinite variance: indeed, if we denote by $L$ the length of a
rise, 
$$\begin{array}{rl}P^{\mu^H}(L\geq n)&=\int_0^1P^p(L\geq
n)d\mu^H(p)=\int_0^1p^nd\mu^H(p)\\&=\ds(2-2H)\int_0^1
p^n(1-p)^{1-2H}dp\ds\mathop{\sim}_{n\to\infty}\Gamma(3-2H){1\over n^{2-2H}}
\end{array}$$ 

(for $H=1/2$, we get $P^{\mu^{1\over2}}(L\geq
n)\ds\mathop{\sim}_{n\to\infty}{1\over n}$)

and we get the same kind of tail as in \cite{TWS}.

The difference lies in the fact that the falls and the rises are not independent: indeed, if
the first rise is short, it probably means that the environment $p$ is small, so that the
following fall will be probably short also. More precisely, if $L_1$ and $L_2$ denote the
lengths of the first rise (resp. fall) and of its following fall (resp. rise),
$$\begin{array}{rl}P^{\mu^H}(L_1\geq n,L_2\geq m)&=E_{\mu^H}[P^p(L_1\geq
n,L_2\geq
m)]=E_{\mu^H}[p^{n+m}]\\&\ds\mathop{\sim}_{n\to\infty}\Gamma(3-2H){1\over
(n+m)^{2-2H}}
\end{array}$$
$L_1$ and $L_2$  are therefore not independent.

Note: Theorems 1 and 2 can be extended to any probability measure with moments
equivalent to ${1\over n^{2-2H}L(n)}$, where $L$ is a slowly varying function. The
additional arguments can be found in \cite{Taq} and rely mainly on Karamata's theorem
in order to replace the naive comparison between sums and integrals we used at the end of
the proof. 
We chose $\mu^H$ in our statement for the simple formula it yields for $c_H$ and because it
is the law of  $1-{U^{1\over2-2H}\over2}$ where $U$ is uniform on $[0,1]$, which makes
it easy to simulate. At the end of the article, we will discuss the practical interest of
other measures.

\section{The case $0<H<1/2$}
\subsection{The alternating correlated random walk}
We first remark that the correlated random walks of section 2, cannot
provide negative correlations for the increments, at least for increments separated by an
even time interval. The best we can hope is to get an alternate sign for the correlations.
In order to get a process with always negative correlations (except for variances), we will
consider the sequence of the sum of two consecutive increments.  
More precisely, if we consider, with the notations of section 2, the sequence 
$(\eps^p_{2n+1}+\eps^p_{2n+2})_{n\geq0}$, for any $p$ less than $1/2$, we get indeed a
sequence of negatively correlated variables and it is also possible to exhibit a probability
on $[0,1]$ such that the equivalence relation (1) mentioned in the introduction will be
satisfied.

But, in the case $0<H<1/2$, this condition alone does not ensure a
scaling limit for the scale $N^H$ (which is that time smaller than $\sqrt N$). 
It has to be allied to a compensation relation between all the correlations. We refer to
\cite{Taq} (section 5) for the statement of this condition, and we will explicit it
further.  

It is the reason why we have to introduce a different kind
of walk, we will call alternating correlated random walk with persistence $p$. It is a
process evolving on $\Z$ by jumps of +1 or -1, whose probability of making the same jump as
the previous one is alternately $p$ and 0. In other words, one jump over two is the
opposite of the previous one:

\begin{defi}
For any $p\in[0,1]$, the alternating correlated random walk $\t X^p$ with persistence $p$,
is a
$\Z$-valued discrete process, such that :

\_ $\t X_0^p=0$, $P(\t X_1^p=-1)=1/2$, $P(\t X_1^p=1)=1/2$. 

\_ $\forall n\geq1, \, \t\eps_n^p:=\t X_{n}^p-\t X_{n-1}^p$ equals 1 or -1 a.s.

\_ $\forall n\geq1, P(\t\eps_{2n}^p=\t\eps_{2n-1}^p|\sigma(\t X_k^p, 0\leq k\leq
2n-1))=p.$

\_ $\forall n\geq1, \t\eps_{2n+1}^p=-\t\eps_{2n}^p$.

\end{defi}

As suggested in the introduction of this section, we will be actually interested by the
process $$Y_n^p:={\t X_{2n}^p\over2\sqrt{p}}$$ (the importance of this normalization will
appear later). The trajectories of this  process take only two values, which are either
$-1/\sqrt p$ and 0 or 0 and $1/\sqrt p$. The successive lengths of the time intervals
during which the process stays on each value are independent and  geometric with parameter
$1-p$. 

We compute now the correlations of the increments of $Y_n^p$ i.e. of $\delta_n^p:=
Y_{n}^p-Y_{n-1}^p={1\over 2\sqrt p}(\t\eps_{2n-1}^p+\t\eps_{2n}^p)$ for $n\geq1$.

\begin{prop}\_ $\forall m\geq1,$
$E[(\d_{m}^p)^2]=1  $

\_ $\forall m\geq1, n\geq1,$
$E[\d_{m}^p\d_{m+n}^p]=-p(1-2p)^{n-1}  $

\end{prop}
Proof:
As in Proposition 1, everything is based on the following facts:

 $\forall n\geq1$,

\_ $E[\t\eps_{2n}^p|\sigma(\t X_k^p, 0\leq k \leq 2n-1)]=(2p-1)\t\eps_{2n-1}^p$,

\_ $E[\t\eps_{2n+1}^p|\sigma(\t X_k^p, 0\leq k \leq 2n)]=-\t\eps_{2n}^p$,

$\forall m\geq1,$
\par\noindent
$\begin{array}{rl}E[(\d_{m}^p)^2]&={1\over
4p}(E[(\t\eps_{2m-1}^p)^2]+E[(\t\eps_{2m}^p)^2]+2E[\t\eps_{2m-1}^p\t\eps_{2m}^p])\\
&={1\over4p}(2+2(2p-1))\\
&=1\end{array}$

$\forall m\geq1, n\geq1,$ by successive conditionings,
\par\noindent
$\begin{array}{rl}E[\d_{m}^p\d_{m+n}^p]&={1\over
4p}E[(\t\eps_{2m-1}^p+\t\eps_{2m}^p)(\t\eps_{2m+2n-1}^p+\t\eps_{2m+2n}^p)]\\
&={1\over4p}(1-2p)^{n-1}
E[(\t\eps_{2m-1}^p+\t\eps_{2m}^p)(\t\eps_{2m+1}^p+\t\eps_{2m+2}^p)]\\
&={1\over4p}(1-2p)^{n-1}(-1-2(2p-1)-(2p-1)^2)\\&=-p(1-2p)^{n-1} \end{array}$

 \qed

We can already note the following (compensation) relation: 
$$E[(\d_{m}^p)^2]+2\sum_{n\geq1}E[\d_{m}^p\d_{m+n}^p]=0\qquad(2)$$

Again we introduce an extra randomness in the persistence.
 We first denote by $Q^p$ the law of $Y^p$
for a given $p$. Now, considering a probability measure $\mu$ on $[0,1]$, we will call
$Q^\mu$, the annealed law of the correlated walk associated to
$\mu$, i.e. the measure on $\Z^\N$ defined by $dQ^\mu:=\int_0^1 Q^p d\mu(p)$. 

Let $Y^\mu$ be a process of law $Q^\mu$.
We introduce  $\d^\mu_n:=
Y_{n}^\mu-Y_{n-1}^\mu$. From Proposition 3 we get the straightforward result: 
\begin{prop}
\_ $\forall m\geq1,$
$E[(\d_{m}^\mu)^2]=1. $

\_ $\forall m\geq1, n\geq1,$
$E[\d_{m}^\mu\d_{m+n}^\mu]=-\int_0^1p(1-2p)^{n-1}d\mu(p). $
\end{prop}
Note: As the compensation relation (2) is satisfied "$p$ by $p$", it remains true for the
annealed correlations.

The goal now is to introduce a probability measure $\mu$ leading to the same
equivalent as (1), mentioned in the introduction.

\subsection{Statement and proof of the result}
We  proceed as in previous section:
\begin{theo}
Let $H\in]0,1/2[$. 

Denote by $\mu^H$ the probability on $[0,{1\over2}]$ with density
$(1-2H)2^{1-2H}p^{-2H}$.

Let $(Y^{\mu^H,i})_{i\geq1}$  be a sequence of independent processes of law $P^{\mu^H}$,
$$ \ml^\md\lim_{N\to\infty}\ml\lim_{M\to\infty} c_H
{Y^{\mu^H,1}_{[Nt]}+...+Y^{\mu^H,M}_{[Nt]}\over N^H\sqrt{M}} = B_H(t)$$
with $c_H=\sqrt{2H\over\Gamma(2-2H)}$.
\end{theo}
Proof: 
The central limit theorem implies that $ \ml\lim_{M\to\infty} 
{Y^{\mu^H,1}_{k}+...+Y^{\mu^H,M}_{k}\over \sqrt{M}}$ is a discrete centered Gaussian process
$(Z^H_k)_{k\geq1}$, with stationary increments $W^H_k:=Z^H_{k+1}-Z^H_k$ with $E[W^H_k]=0$,
$E[(W^H_k)^2]=1$ and 
$$\forall i,n\geq0,\quad r(n):=E[W^H_i W^H_{i+n}]=-(1-2H)2^{1-2H}\int_0^{1\over2}(1-2u)^n
u^{1-2H}du$$ 
$\begin{array}{rl}
r(n)&=-\ds{(1-2H)\over2}\int_0^1(1-v)^n v^{1-2H}dv\\
&=-\ds{(1-2H)\over2}{\Gamma(n+1)\Gamma(2-2H)\over\Gamma(n+3-2H)}\\
&\ds\mathop{\sim}_{n\to\infty} -{(1-2H)\over2}\Gamma(2-2H){1\over n^{2-2H}}
={1\over c_H^2}{H(2H-1)\over n^{2-2H}}\end{array}$

So that,
 
$\begin{array}{rl}E[c_H^2(G^H_1+...+G^H_N)^2]&=\ds c_H^2\sum_{i=1}^N\sum_{j=1}^N r(|i-j|)\\
&=
\ds c_H^2(r(0)+\sum_{i=1}^{N-1}[r(0)+2\sum_{k=1}^i r(k)])\\
&= \ds c_H^2(r(0)-2\sum_{i=1}^{N-1}\sum_{k=i+1}^\infty r(k))\\
&\ds\mathop{\sim}_{n\to\infty} N^{2H}\end{array}$

The last equality comes from the compensation relation (2), and 
the last step consists simply in two successive comparisons between sums and integrals.
A direct application of \cite{Taq} (lemma 5.1) allows to conclude.
\qed

The process $Y_n^p$ reminds strikingly the On/Off processes of the traffic  modeling
theory, described in \cite {TWS}, or the $"V"$-renewal process of \cite{TL} which would have
alternate rewards. But our attempt to make this model work here, failed probably because
it is not clear how to translate to this model the normalization of $\t X_{2n}^p$ by
$\sqrt{p}$, which is crucial in our construction because its absence would yield an
infinite measure for
$\mu^H$.

Remark 1: As in previous section, $\mu^H$ is not the only convenient measure. We remark that
it is the law of ${U^{1\over1-2H}\over2}$, where $U$ is uniform on $[0,1]$.

Remark 2: Other normalizations of $\t X_{2n}^p$ by powers of $p$ can be considered, but then
the measure $\mu^H$ has to be multiplied by a power of $p$ also.

\section{Practical aspects}
As mentioned above, taking the limit in reverse order, yields trivial processes. The
question is: for a given number of steps $N$, what kind of $M$ are we supposed to take in
order to approximate the right limit ? We restrict our study to the law of $B_H(1)$. We
base our study on  Berry-Esseen's inequality,
applied to the sequence of i.i.d. variables with law
$c_HX_N^{\mu^H}/N^H$ (resp.
$c_HY_N^{\mu^H}/N^H$) for
$H\geq1/2$ (resp. $H<1/2$).

What remains to compute, is an upper bound for the third moment
of the absolute value of these variables.

We begin with the case $H>{1\over2}$. 

\begin{prop}
For $H>1/2$, for $N$ large  enough,
$$E[({c_H|X_N^{\mu^H}|\over N^H})^3]\leq D_H
N^{1-H},$$

with $D_H=\sqrt{6(2H-1)\over(H+1)(2H+1)}\times c_H.$

\end{prop}
Nota Bene: we express $D_H$ in terms of $c_H$, in order to have a formula that works for
other measures, as it will be useful in the following.

Proof: We omit here the superscipts $\mu^H$  in the variables
$X_n^{\mu^H}$'s and $\eps_n^{\mu^H}$'s.
$X_N=\sum_{k=1}^N\eps_k$,  with:

\_ $\eps_n$ are Bernoulli(1/2),

\_ $c_H^2 Cov(\eps_k,\eps_l)=r(|k-l|)$, with $ r(n)\sim_{\infty}{H(2H-1)\over
n^{2-2H}}$.

We use  Cauchy-Schwarz inequality to get:
$$ \begin{array}{rl}
\ds E[({c_H|X_N|\over N^H})^3]&
\ds\leq E[({c_H|X_N|\over N^H})^2]^{1\over2}E[({c_H|X_N|\over N^H})^4]^{1\over2}\\
&\sim\ds{ c_H^2\over N^{2H}}(\ds\mathop{\sum_{1\leq i_k\leq N}}_{1\leq k\leq 4}
E[\eps_{i_1}\eps_{i_2}\eps_{i_3}\eps_{i_4}])^{1/2}
\end{array} $$
(using that the variance converges to 1 for large $N$).

Assume $i_4\leq i_3\leq i_2\leq i_1$, we get, as in
Proposition 1, by successive conditionings:
$$E[\eps_{i_1}^p\eps_{i_2}^p\eps_{i_3}^p\eps_{i_4}^p]=
(2p-1)^{(i_1-i_2)+(i_3-i_4)}  
$$
So that,
$$E[\eps_{i_1}\eps_{i_2}\eps_{i_3}\eps_{i_4}]= r((i_1-i_2)+(i_3-i_4))  
$$ 
We see by using the equivalence between sums and integrals, as we did at the end of the
proof of Theorem 1, that the sum is equivalent to:
$$ 
4!\times\int_0^N\int_0^{x_1}\int_0^{x_2}\int_0^{x_3}
r((x_1-x_2)+(x_3-x_4))dx_4dx_3dx_2 dx_1
$$

$$
\sim {4!\over c_H^2}\times H(2H-1)\int_0^N\int_0^{x_1}\int_0^{x_2}\int_0^{x_3}
((x_1-x_2)+(x_3-x_4))^{2H-2}dx_4dx_3dx_2 dx_1
$$

$$
={6\over c_H^2}{(2H-1)N^{2H+2}\over(H+1)(2H+1)}
$$
 \qed

Now applying Berry-Esseen's inequality, the error of the distribution function of the
marginal at time 1, is dominated by $ 0.65\times D_H{N^{1-H}\over\sqrt{M}} $
(using the constant 0.65 of \cite{Z}, as far as the third moment is much bigger than the
power ${3\over2}$ of the variance). 
 
We deduce that  Theorem 1 remains true as soon as
$M(N)\to\infty$ as $N\to\infty$ at a faster rate than $N^{2-2H}$:

\begin{cor}
Let $M$ be a function on the integers such that ${M(N)\over N^{2-2H}}$ tends to $\infty$, 
$$ \ml^\md\lim_{N\to\infty} c_H
{X^{\mu^H,1}_{[Nt]}+...+X^{\mu^H,M(N)}_{[Nt]}\over N^H\sqrt{M(N)}} = B_H(t)$$

\end{cor}
Proof: Using the generalization of Berry Esseen's inequality to multidimensional variables
\cite{BR}\cite{Saz}, we get the convergence for finite-dimensional marginals. 

To get the weak convergence, we cannot use directly \cite{Taq} as in previous section, and
we prove the tightness of the family of processes, by checking Billingsley's criteria
\cite{Bil} (Theorem 15.6):

Denote by $S_N(t):= {X^{\mu^H,1}_{[Nt]}+...+X^{\mu^H,M(N)}_{[Nt]}\over \sqrt{M(N)}}$.

Let $1\geq t_2\geq t\geq t_1\geq0$, and $k\in\N$,

$\begin{array}{rl}
J_N(k,t_2,t,t_1)&:=E[|{S_N(t_2)-S_N(t)\over N^H}|^k|{S_N(t)-S_N(t_1)\over N^H}|^k]\\
&\leq {1\over N^{2kH}}E[S_N(t_2-t)^{2k}]^{1\over2}E[S_N(t-t_1)^{2k}]^{1\over2}
\end{array}$

$$\forall k\in\N,\quad  E[S_N(t)^{2k}]={1\over
M(N)^k}E[(X^{\mu^H,1}_{[Nt]}+...+X^{\mu^H,M(N)}_{[Nt]})^{2k}]$$

When we develope the polynomial with degree $2k$ inside the expectation, we notice that
only the monomials of the type $(X^1)^{2\alpha_1}...(X^{M})^{2\alpha_M}$ have a non null
contribution, as far as the $X^i$'s are centered and independent. This set has cardinal
$O(M^k)$. 

Now, we find, by the same means than in Proposition 1, that
$$\forall\alpha_1\geq1,\quad
E[X_N^{2\alpha_1}]=O(N^{2\alpha_1+(2H-2)}).$$
This leads to $  E[S_N(t)^{2k}]=O( (Nt)^{2k+(2H-2)k})=
O((Nt)^{2Hk}).  $

Hence, for some positive constant $C$,

 $J_N(k,t_2,t,t_1)\leq
C(t_2-t)^{2Hk}(t-t_1)^{2Hk}\leq C(t_2-t_1)^{2Hk}$.

We choose $k>{1\over 2H}$ in order to satisfy Billingsley's criteria.\qed

We conjecture that
this result remains true in the case of \cite{TWS} which proposes an answer to the question
asked at the end of \cite{TWS}. 

We do the same for the case $H={1\over2}$:

\begin{prop}
For $N$ large  enough,
$$E[( {c_{1\over2}|X_N^{\mu^{1\over2}}|\over \sqrt{N\log N}})^3]\leq
 c_{1\over2}\times{\sqrt{2N}\over\log N}.$$

\end{prop}

Proof: As in previous proposition, 
$X_N=\sum_{k=1}^N\eps_k$,  with:

\_ $\eps_n$ are Bernoulli(1/2),

\_ ${c^2_{1\over2}} Cov(\eps_k,\eps_l)=r(|k-l|)$, with $ r(n)\sim_{\infty}{1\over
2n}$.

Similarly as in previous proposition,
$\ds E[({c_{1\over2}|X_N|\over \sqrt{N\log N}})^4]^{1\over2}$ is equivalent to:
$$\sim 
{c^2_{1\over2}\over N\log
N}({4!\over c^2_{1\over2}}\int_0^N\int_0^{x_1}\int_0^{x_2}\int_0^{x_3}
r((x_1-x_2)+(x_3-x_4))dx_4dx_3dx_2 dx_1)^{1\over2}
$$
But,$$
\int_0^N\int_0^{x_1}\int_0^{x_2}\int_0^{x_3}
{1\over(x_1-x_2)+(x_3-x_4)}dx_4dx_3dx_2 dx_1={N^3\over6}.
$$
\qed

Now applying Berry-Esseen's inequality, the error of the distribution function of the
marginal at time 1, is dominated by $ 1.3\times {\sqrt{N}\over\log N \sqrt{M}} $.

\begin{cor}
Let $M$ be a function on the integers such that ${M(N)\over N\log(N)^{-2}}$ tends to
$\infty$, 
$$ \ml^\md\lim_{N\to\infty} c_{1\over2}
{X^{\mu^{1\over2},1}_{[Nt]}+...+X^{\mu^{1\over2},M(N)}_{[Nt]}\over \sqrt{M(N)\times
N\log(N)}} = B(t)$$

\end{cor}

Finally, we treat the case $H<{1\over2}$:
\begin{prop}
For $H<1/2$, for $N$ large  enough,
$$E[({c_H|Y_N^{\mu^H}|\over N^H})^3]\leq D_H
N^{{1\over2}-H},$$

with $D_H=\sqrt{2H\over2H+1}\times c_H.$

\end{prop}

Proof:
The proof is  more delicate than in Proposition 5, because the computation of the
fourth moment  makes appear compensations (as in the variance computation), and we have
to treat them first "$p$ by $p$", before integrating against $\mu^H$.

Assume $i_4\leq i_3\leq i_2\leq i_1$, we get, as in
Proposition 1, by successive conditionings:
$$E[\d_{i_1}^p\d_{i_2}^p\d_{i_3}^p\d_{i_4}^p]=
r_p(i_1-i_2)r_p(i_3-i_4), $$
where $r_p(0)=1$ and $\forall n>0, \, r_p(n)=-p(1-2p)^{n-1}$.

Now,
$$\begin{array}{rl}
\ds\mathop{\sum_{1\leq i_k\leq N}}_{1\leq k\leq
4}E[\d_{i_1}^p\d_{i_2}^p\d_{i_3}^p\d_{i_4}^p]
&=12\ds\sum_{\max\{i_3,i_4\}\leq\min\{i_1,i_2\}}
E[\d_{i_1}^p\d_{i_2}^p\d_{i_3}^p\d_{i_4}^p]+O(N)\end{array}$$
$$
=12\ds\sum_{\max\{i_3,i_4\}\leq\min\{i_1,i_2\}}r_p(|i_1-i_2|)r_p(|i_3-i_4|)+O(N)$$
$$
=12\ds\sum_{\min\{i_1,i_2\}=1}^N\left(\sum_{\max\{i_3,i_4\}=1}^{\min\{i_1,i_2\}}
r_p(|i_3-i_4|)\right)r_p(|i_1-i_2|)+O(N)$$
 $$
=12\ds\sum_{\min\{i_1,i_2\}=1}^N\left(\sum_{\max\{i_3,i_4\}=1}^{\min\{i_1,i_2\}}
\left(r_p(0)+2\sum_{l=1}^{\max\{i_3,i_4\}}r_p(l)\right)\right)r_p(|i_1-i_2|)+O(N)$$
$$
=12\ds\sum_{\min\{i_1,i_2\}=1}^N\left(\sum_{\max\{i_3,i_4\}=1}^{\min\{i_1,i_2\}}
(1-2p)^{\max\{i_3,i_4\}}\right)r_p(|i_1-i_2|)+O(N)
$$
$$=12\ds\sum_{\min\{i_1,i_2\}=1}^N\left(\sum_{k=1}^{\min\{i_1,i_2\}}
(1-2p)^{k}\right)\left(r_p(0)+2\sum_{l=1}^{N-\min\{i_1,i_2\}}r_p(l)\right)+O(N)$$
$$=12\ds\sum_{\min\{i_1,i_2\}=1}^N\left(\sum_{k=1}^{\min\{i_1,i_2\}}
(1-2p)^{k}\right)(1-2p)^{N-\min\{i_1,i_2\}}+O(N)$$
$$=12\ds\sum_{1\leq k\leq l\leq N}(1-2p)^{N-(l-k)}+O(N)
$$
(The contribution of the exceptional situations $\max\{i_3,i_4\}=\min\{i_1,i_2\}$ is
estimated by $O(N)$, because of the same compensations as the ones described in the above
equalities.)

 We now use,
 $$\begin{array}{rl}
\int_0^1(1-2p)^nd\mu^H(p)&\ds=-2\sum_{k>n}r(k)
\\&\ds\sim{1\over c_H^2}{2H\over n^{1-2H}}
\end{array}$$

We deduce,

$$ {c_H^2\over N^{2H}}(\ds\mathop{\sum_{1\leq i_k\leq N}}_{1\leq k\leq
4}E[\d_{i_1}\d_{i_2}\d_{i_3}\d_{i_4}])^{1\over2} \sim$$
$$
{c_H^2\over N^{2H}}
\left({2H\over c_H^2}\int_0^N\int_0^{x_2}{dx_1\over(N-(x_2-x_1))^{1-2H}}\right)
^{1\over2}
$$
$$=c_H\times\sqrt{2H\over2H+1}N^{{1\over2}-H}.$$
\qed

Applying Berry-Esseen's inequality, the error of the distribution function of the
marginal at time 1, is dominated by $ 0.65\times D_H {N^{{1\over2}-H}\over\sqrt{M}} $.

\begin{cor}
Let $M$ be a function on the integers such that ${M(N)\over N^{1-2H}}$ tends to $\infty$, 
$$ \ml^\md\lim_{N\to\infty} c_H
{X^{\mu^H,1}_{[Nt]}+...+X^{\mu^H,M(N)}_{[Nt]}\over N^H\sqrt{M(N)}} = B_H(t)$$
\end{cor}

As a general conclusion, we can say that, using $\mu_H$, the number of computations we have
to make, in order to get a $N$ steps trajectory, is a constant times $N^{3-2H}$ for
${1\over2}<H$, and 
$N^{2-2H}$ for $0<H<{1\over2}$. In any case, it is a power of $N$ between 1 and 2.

The constant in factor, is of big importance for practical simulation. We took the best
constant in the Berry-Esseen's inequality we could find in the literature, even if it is 
bigger than the best constant possible conjectured by Esseen, i.e.
${3+\sqrt{10}\over6\sqrt{2\pi}}\simeq0.41$ that would gain in time a squared factor equal
to 2.5 (see \cite{Z} for a nice discussion on this subject). 

But more important is to note that the constant can be considerably ameliorated by using
other measures than $\mu_H$, providing smaller $c_H$'s. 
It is the case for the real-indexed sequence of probabilities $(\mu_{H,k})_{k>0}$, defined
below:

\_  $\forall
H\in[{1\over2},1],$
$$d\mu_{H,k}(p):=2^{k+1-2H}{\Gamma(k+2-2H)\over\Gamma(k)\Gamma(2-2H)}
(p-{1\over2})^{k-1}(1-p)^{1-2H}1_{[{1\over2},1]}(p)dp$$ which is just the law of
${1+B(k,2-2H)\over2}$ and coincides with $\mu_H$ for $k=1$ (where $B(a,b)$ denotes the
Beta variable with parameters $a$ and $b$). 

\_ $\forall
H\in]0,{1\over2}[,$
$$d\mu_{H,k}(p):=2^{k-2H}{\Gamma(k+1-2H)\over\Gamma(k)\Gamma(1-2H)}
({1\over2}-p)^{k-1}p^{-2H}1_{[0,{1\over2}]}(p)dp$$ which is just the law of
${B(1-2H,k)\over2}$ and coincides with $\mu_H$ for $k=1$.

\_ For
$H>{1\over2}$, we obtain
$c_{H,k}=\sqrt{H(2H-1)\Gamma(k)\over\Gamma(k+2-2H)}\sim {\sqrt{H(2H-1)}\over k^{1-H}}$.

So that, using $\mu_{H,k}$ for large $k$, yields an error equivalent to: 
$$0.65\sqrt{6H(2H-1)^2\over(H+1)(2H+1)}\times({N\over k})^{1-H}{1\over\sqrt{M}}$$

\_ For
$H={1\over2}$, we obtain
$c_{{1\over2},k}={1\over\sqrt{2k}}$.

The error is equivalent to: 

$${0.65 \over\log N}\times({N\over k})^{1\over2}{1\over\sqrt{M}} $$

\_ For
$H<{1\over2}$, we obtain
$c_{H,k}=\sqrt{2H\Gamma(k)\over\Gamma(k+1-2H)}\sim {\sqrt{2H}\over k^{{1\over2}-H}}$.

The error is equivalent to: 

$$0.65 \sqrt{4H^2\over2H+1}\times({N\over k})^{{1\over2}-H}{1\over\sqrt{M}} $$

Using $\mu_{H,k}$ instead of $\mu_H$, allows to gain for $M$ a factor $k^{2-2H}$ (resp.
$k^{1-2H}$) for $H>{1\over2}$ (resp. for $H<{1\over2}$). Loosely speaking, it erases the
noise generated by the parameters between ${1\over2}$ and a fixed constant smaller than 1. 
The trouble making
$k$ increase, is that it damages the value of the covariance of $X_N$ (resp. $Y_N$), but 
we can allow any $k(N)=o(N)$. So that we obtain an algorithm with any $M(N)=1/o(1)$ number
of trajectories. As a result, our algorithm requires a number of computations of the order
$N/o(1)$, for any $o(1)$. Moreover, we have only $M(N)$ real datas to keep in
memory along the whole procedure, corresponding to the coin-tossed parameters of the walks,
and $M(N)$ integers (-1 or 1), giving the last moves of the walks.

We want to give now a second family of measures $(\mu'_{H,k})_{k>0}$:

\_ For $H>{1\over2}$, $\mu'_{H,k}$ is the law of $1-{(1-U^{1\over k})^{1\over 2-2H}\over2}$.

An easy computation gives $c'_{H,k}={c_H\over\sqrt{k}}$ (error: $0.65\times
D_H{N^{1-H}\over\sqrt{kM}}$).

\_ For $H<{1\over2}$, $\mu'_{H,k}$ is the law of ${(1-U^{1\over k})^{1\over 1-2H}\over2}$.

Again, an easy computation gives $c'_{H,k}={c_H\over\sqrt{k}}$ (error: $0.65\times
D_H{N^{{1\over2}-H}\over\sqrt{kM}}$).

The advantage of this family is obviously the easy simulation it provides. The error is
estimated by a term containing ${1\over\sqrt{k}}$, that seems to be better than the last
one, but the damages on the variances grow faster than for $\mu_{H,k}$. Actually the scale
$\sqrt{k}$ corresponds, in the previous family, to the scale $k^{1-H}$ (resp.
$k^{{1\over2}-H}$) for $H>{1\over2}$ (resp. $H<{1\over2}$).
The drawback of this family, is that the theoretical computations of the variance are not
very explicit.

   We illustrate our results by three graphs
corresponding to  three different parameters of $H$, with $N=1000$, and a theoretical error
smaller than $10\%$ (we indicate the time it takes for Matlab to draw a graph):

\_ For $H=0.25$,  we take $M=200$ and $k=1$ (15 seconds)  (Fig. 1).

\_ For $H=0.5$, we can see that the
convergence in Theorem 2 is the slowest one, and we will use the simple random walk to
simulate it. (Fig. 2)  

\_ For $H=0.75$,  we take $M=400$ and $k=0.5$ (25 seconds) (Fig. 3). 

Remark: The interest of taking large $k$ appears when $N$ is very large, especially if we
want to preserve a small error on the whole trajectory. In the case $H>{1\over2}$, we were
even  obliged to take $k$ smaller than one.

We remark the different behaviours of
the trajectories: let us recall that the Hausdorff dimension of the trajectories are a.s.
equal to $2-H$, and we notice that the variances of the process between 0 and 1 become
larger when $H$ decreases.


As we noticed just above, we may be limited by the fidelity of the covariance of our
process. In this spirit, it is quite interesting to note that the  autocovariance function
of the Gaussian noise $G_H(j)$ of the introduction is, up to a shift,
 the sequence of  moments of a probability measure on $[0,1]$. I first remarked it by
checking the conditions of Hausdorff theorem (\cite{F} p.226), but Marc Yor gave me kindly,
the method to get the explicit density of this measure, and I present it here. This brings
a third (the last !) family of probability measures:

\begin{prop}
Let $H\in]{1\over2},1[$.
Consider the family of probability measures $(\nu_{H,k})_{k>0}$ on $[{1\over2},1]$, with
density
$C(H,k)\times (1-p)^{2}(2p-1)^{k-1}(\ln({1\over2p-1}))^{-1-2H}$,
and 

$C(H,k):={16H(2H-1)\over\Gamma(2-2H)}\times ((k+2)^{2H}-2(k+1)^{2H}+k^{2H})^{-1}  $
$$\forall n\geq0, \int_0^1 (2p-1)^n
d\nu_{H,k}(p)={(n+k+2)^{2H}-2(n+k+1)^{2H}+(n+k)^{2H}\over (k+2)^{2H}-2(k+1)^{2H}+k^{2H}}  
$$

\end{prop} 

Proof: $\forall n\geq 1$,
$$
{1\over2}((n+1)^{2H}-2n^{2H}+(n-1)^{2H})=
\ds H\int_0^1(n+t)^{2H-1}-(n+t-1)^{2H-1} dt$$
$$
=H(2H-1)\int_0^1\int_0^1(n+t+s-1)^{2H-2}dsdt$$
$$={H(2H-1)\over\Gamma(2-2H)}\int_0^1\int_0^1\int_0^{+\infty}e^{-(n+t+s-1)u}u^{1-2H}dudsdt$$

$$ ={H(2H-1)\over\Gamma(2-2H)}\int_0^{+\infty}e^{-nu}(1-e^{-u})^2e^uu^{-1-2H}du$$
$$={H(2H-1)\over\Gamma(2-2H)}\int_0^1x^{n}({1-x\over x})^2(\ln{1\over x})^{-1-2H}dx  $$
 The end of the proof is straightforward, by change of variable.\qed

Note: The normalization corresponding to  $\nu_{H,k}$ is:

$c''_{H,k}=({(k+2)^{2H}-2(k+1)^{2H}+k^{2H}\over2})^{1\over2}\sim {\sqrt{H(2H-1)}\over
k^{1-H}}\sim c_{H,k} $.

The same can be done for $H\in]0,{1\over2}[$:
\begin{prop}
Let $H\in]0,{1\over2}[$.
Consider the family of probability measures $(\nu_{H,k})_{k>0}$ on $[0,{1\over2}]$, with
density
$C(H,k)\times p(2p-1)^{k-1}(\ln({1\over1-2p}))^{-1-2H}$,
and 

$C(H,k):={8H\over\Gamma(1-2H)}\times ((k+1)^{2H}-k^{2H})^{-1}  $
$$\forall n\geq1, -\int_0^1 p(2p-1)^{n-1}
d\nu_{H,k}(p)={(n+k+1)^{2H}-2(n+k)^{2H}+(n+k-1)^{2H}\over2
((k+1)^{2H}-k^{2H})}   $$

\end{prop}
Proof: $\forall n\geq 1$,

$$
{1\over2}((n+1)^{2H}-2n^{2H}+(n-1)^{2H})=
\ds H\int_0^1(n+t)^{2H-1}-(n+t-1)^{2H-1} dt$$

$$={H\over\Gamma(1-2H)}\int_0^1\int_0^{+\infty}
(e^{-(n+t)u}-e^{-(n+t-1)u})u^{-2H}dudt$$

$$ =-{H\over\Gamma(1-2H)}\int_0^{+\infty}e^{-(n+1)u}(e^{u}-1)^2u^{-1-2H}du$$
$$=-{H\over\Gamma(1-2H)}\int_0^1x^{n}({1-x\over x})^2(\ln{1\over x})^{-1-2H}dx  $$
$$=-{8H\over\Gamma(1-2H)}\int_0^{1\over2}(1-2p)^{n-2}p^2(\ln{1\over
1-2p})^{-1-2H}dp 
$$
We find the normalization, using
$\sum_{n\geq1}(n+k+1)^{2H}-2(n+k)^{2H}+(n+k-1)^{2H}=k^{2H}-(k+1)^{2H}$, and the relation
$1=2\sum_{n\geq1}p(1-2p)^{n-1}$.\qed

Note: The normalization corresponding to  $\nu_{H,k}$ is:

$c''_{H,k}=((k+1)^{2H}-k^{2H})^{1\over2}\sim {\sqrt{2H}\over k^{{1\over2}-H}}\sim
c_{H,k}
$.

These  measures $\nu_{H,k}$ have an interest, as far as $X_N^{\nu_H}$ (resp.
$Y_N^{\nu_H}$) have explicit covariances matching quite well with the covariances of the
fractional Brownian motion. On the other hand, they don't seem very easy to simulate. 

\centerline{\epsfig{file=fig1.eps, scale =0.36}}

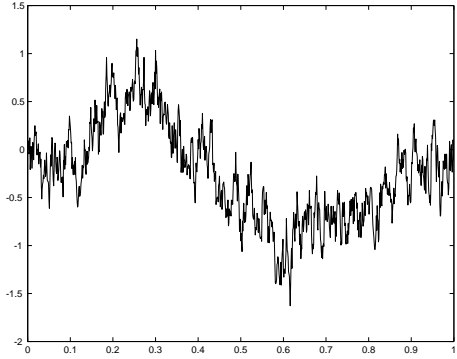
\centerline{\small Fig. 1: H=0.25}

\centerline{\epsfig{file=fig2.eps, scale =0.36}}

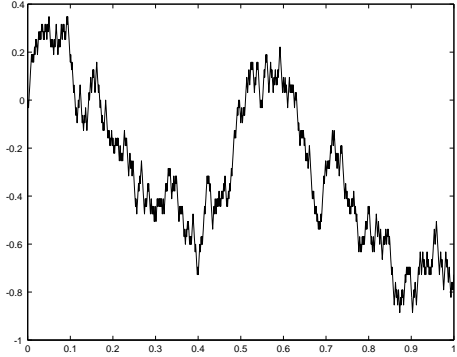
\centerline{\small Fig. 2: H=0.5}

\centerline{\epsfig{file=fig3.eps, scale =0.36}}

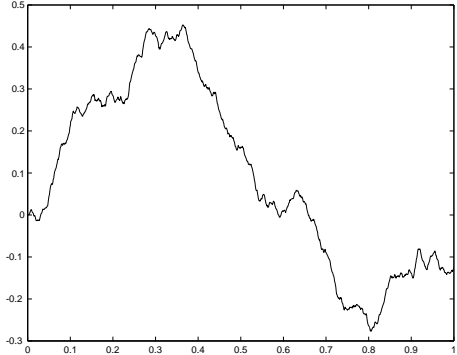
\centerline{\small Fig. 3: H=0.75}

{\bf Aknowledgements:} First of all, I would like to thank Erick Herbin for asking me the
question about the relation between fractional Brownian motion and correlated random walks.
It is a pleasure for me to thank also Marc Yor who helped me finding $\nu_H$, and Zhan Shi
for fruitful conversations.

\end{document}